\documentclass[graybox, envcountchap]{svmult}

\usepackage{makeidx}         
\usepackage{graphicx}        
\usepackage{multicol}        
\usepackage[bottom]{footmisc}

\usepackage{newtxtext}        
\usepackage{newtxmath}       
\usepackage{amsmath}
\usepackage{hyperref}

\makeindex             

\def \calT {T}
\def \calO {\mathcal{O}}
\def \calM {\mathcal{M}}
\def \calH {\mathcal{H}}
\def \calV {\mathcal{V}}
\def \R {\mathbb{R}}
\newcommand{\psd}[2]{\mathcal{S}_+(#1, #2)}

\def \M {\mathcal{M}}
\def \bspline {\mathbf{B}}

\begin{document}

\frontmatter



\mainmatter
%
%
%

%
%
%
%
%
%
%

\title{Balanced truncation for parametric linear systems
	using interpolation of Gramians: a~comparison of algebraic and geometric approaches}
\titlerunning{Parametric balanced truncation by interpolation of Gramians} 
\author{N.T. Son, P.-Y. Gousenbourger, E. Massart, and T. Stykel}
\institute{N. T. Son \at ICTEAM, UCLouvain, Avenue Georges Lema\^itre 4-6/L4.05.01, 1348 Louvain-la-Neuve, Belgium and Thai Nguyen University of Sciences, 25000 Thai Nguyen, Vietnam. \email{thanh.son.nguyen@uclouvain.be}
\and P.-Y. Gousenbourger \at ICTEAM, UCLouvain, Avenue Georges Lema\^itre 4-6/L4.05.01, 1348 Louvain-la-Neuve, Belgium. \email{pierre-yves.gousenbourger@uclouvain.be}
\and E. Massart \at NPL-postdoctoral research associate, Mathematical Institute, University of Oxford, OX2 6GG, UK. \email{estelle.massart@maths.ox.ac.uk}. Most of this work was done when this author was with ICTEAM, UCLouvain. 
\and T. Stykel \at Institute of Mathematics, 
University of Augsburg, 
Universit\"atsstr.~14,
D-86159 Augsburg, Germany. \email{stykel@math.uni-augsburg.de}}
%
%
\maketitle

\abstract{When balanced truncation is used for model order reduction, one has to solve a pair of Lyapunov equations for two Gramians and uses them to construct a reduced-order model. Although advances in solving such equations have been made, it is still the most expensive step of this reduction method. Parametric model order reduction aims to determine 
reduced-order models for parameter-dependent systems. Popular techniques for parametric model order reduction rely on interpolation. Nevertheless, interpolation of Gramians is rarely mentioned, most probably due to the fact that Gramians are symmetric positive semidefinite 
 matrices, a property that should be preserved by the interpolation method.  In this contribution, we propose and compare two approaches for Gramian interpolation. In the first approach, the interpolated Gramian is 
computed 
as a linear combination of the data Gramians with positive coefficients. Even though 
positive semidefiniteness is guaranteed 
in this method, the rank of the interpolated Gramian can be significantly larger than that of the data Gramians.
The second approach aims to tackle this issue by performing the interpolation on the manifold of fixed-rank positive semidefinite matrices. The results of the interpolation step are then used to construct parametric reduced-order models, which are compared numerically on two benchmark problems.
}

\section{Introduction}\label{Sec:Intro}
The need for increasingly accurate simulations in sciences and
technology results in large-scale mathematical models. Simulation of those systems is usually time-consuming or even infeasible, especially with limited computer resources. Model order reduction (MOR) is a well-known tool to
deal with such problems. Founded about half a century ago, 
 this field is still getting attraction due to the fact that many complicated or large problems have not been considered and many advanced methods have not been invoked yet. 

Often, the full order model (FOM) depends on parameters. The reduced-order model (ROM), preferably parameter-dependent as well, is therefore required to approximate the FOM on a given parameter domain. This problem, so-called parametric MOR (PMOR), has been addressed by various approaches such as Krylov subspace-based \cite{FengRK05,LiBS09}, optimization \cite{BaurBBG11}, interpolation \cite{AmsaF08,BaurB09,PanzMEL10,Son12}, and reduced basis technique \cite{HaasO11,SonS17}, just to name a few. The reader is referred to the survey \cite{BennGW15}  and the contributed book \cite{BennOPRU17} for more details. We focus here on the methods that use interpolation to build a ROM for the linear parametric control system
\begin{equation}
\label{eq:psys}
\arraycolsep=2pt
\begin{array}{rcl}
E(\mu)\dot{x}(t,\mu) & = &A(\mu)x(t,\mu)+B(\mu)u(t),\\
y(t,\mu) & = & C(\mu)x(t,\mu),
\end{array}
\end{equation}
where $E(\mu)$, $A(\mu) \in \mathbb{R}^{n\times n}$, $B(\mu) \in \mathbb{R}^{n\times m}$, 
$C(\mu) \in \mathbb{R}^{p\times n}$ with $p,m \ll n$, and \mbox{$\mu\in\mathcal{D}\subset \mathbb{R}^{\ell}$}. We assume that the matrix 
$E(\mu)$ is nonsingular and all eigenvalues of the pencil $\lambda E(\mu)-A(\mu)$ have negative real part for all $\mu\in\mathcal{D}$. This assumption is to avoid working with singular control systems and to restrict ourselves to the use of standard balanced truncation \cite{Dai89,Styk04}.  The goal is to approximate system \eqref{eq:psys} with a smaller parametric model
\begin{equation}
\label{eq:redpsys}
\arraycolsep=2pt
\begin{array}{rcl}
\tilde{E}(\mu)\dot{\tilde{x}}(t,\mu) & = &\tilde{A}(\mu)\tilde{x}(t,\mu)+\tilde{B}(\mu)u(t),\\
\tilde{y}(t,\mu) & = &\tilde{C}(\mu)\tilde{x}(t,\mu),
\end{array}
\end{equation}
where $\tilde{E}(\mu)$, $\tilde{A}(\mu) \in \mathbb{R}^{r\times r}$, $\tilde{B}(\mu)\in\mathbb{R}^{r\times m}$, 
$\tilde{C}(\mu) \in \mathbb{R}^{p\times r}$ and $r \ll n$. 

Interpolation-based methods work as follows. On a given sample grid $\mu_i$, \linebreak $i = 1, \dots, q$, in the parameter domain $\mathcal{D}$, one computes a ROM associated with each $\mu_i$. These ROMs can be obtained
using any MOR method for non-parametric models \cite{Anto05} and characterized by either their  projection subspaces, coefficient matrices, or transfer functions. Then they are interpolated using standard 
 methods like Lagrange or spline interpolation. These approaches 
have been discussed intensively in many publications see, e.g., \cite{PanzMEL10,DegrVW10,AmsaF11} for interpolating local reduced system matrices, \cite{AmsaF08,Son13} for interpolating projection subspaces,  \cite{BaurB09,SonS15} for interpolating reduced transfer functions, and \cite{Zimm19} for a detailed discussion on the use of manifold interpolation for model reduction.
 Each of them has its own strength and acts well in some specific applications but fails to be superior to the others in a general setting. 

When balanced truncation \cite{Moor81} is used, one has to solve a pair of Lyapunov equations for  two Gramians. Although advances in solving such equations have been made, it is still the most expensive step in this reduction method. Therefore, any interpolation method that can circumvent this step is of interest. Unfortunately, to our knowledge, there has been no work addressing this issue. In this contribution, we propose to interpolate the solutions to these equations, i.e., the Gramians. It is noteworthy that in the large-scale setting, one should avoid working with full-rank solution matrices. Fortunately, in many practical cases, the solution of the Lyapunov equation can be well approximated by symmetric positive semidefinite (SPSD) matrices of considerably smaller rank~\cite{Penz00b,AntoSZ02}. 
Such approximations can be used in the square root balanced truncation method \cite{TombP87}
making the reduction procedure more computationally efficient.

To ensure that the SPSD property is preserved during the interpolation, we propose two approaches. In the first one, which is the main content of Section~\ref{Sec:BT_standard interpolation}, the target Gramians are written as a linear combination of the data Gramians with some given (positive) weights. The main drawback of this approach is that the interpolated Gramians will, in general, have a considerably larger rank than the data Gramians. However, when combining it with the reduction process, we can truncate the unnecessary data and design an offline-online decomposition of the whole procedure, to reduce the computational cost of the operations that have to be done on-the-fly, i.e., that depend on the value of the parameter associated to the target ROM. 
 We refer to this as the linear algebraic (or algebraic for short) approach. 
The second approach, given in Section~\ref{Sec:Manifold}, consists of mapping beforehand all the matrices to the set of fixed-rank 
positive semidefinite matrices, and performing the interpolation directly in that set. This would ensure that the rank of the interpolated Gramians remains consistent with the ranks of the data Gramians. It was shown in \cite{VandAV09,MassA18} that the set of SPSD matrices of fixed rank can be turned into a Riemannian manifold by equipping it with a differential structure. We can then resort to interpolation techniques specifically designed to work on Riemannian manifolds. Oldest techniques are based on subdivision schemes~\cite{Dyn2009} or rolling procedures~\cite{Huper2007}. In the last decades, path fitting techniques rose up, such as least squares smoothing~\cite{Machado2010} or more recently by means of B\'ezier splines~\cite{Absil2016,Gousenbourger2018}. 
The latter will be employed here for interpolating the Gramians. The resulting PMOR
method will be referred to as the geometric method in the sense that it strictly preserves the geometric structure of data. 

The rest of the paper is organized as follows. In Section~\ref{Sec:BT_standard interpolation} we briefly recall balanced  truncation for MOR, the square root balanced truncation procedure, and present the algebraic 
interpolation method. 
Section~\ref{Sec:Manifold} is devoted to the geometric 
interpolation method. It first describes the geometry of the manifold of fixed-rank SPSD matrices, and then algorithms to perform interpolation on this manifold. The two proposed approaches are then compared numerically in Section~\ref{Sec:NumerExam}, and the conclusion is given in Section~\ref{Sec:Concl}. 


\section{Balanced truncation for parametric linear systems and standard interpolation}\label{Sec:BT_standard interpolation}

\subsection{Balanced truncation}
Balanced truncation \cite{Anto05,Moor81} is a well-known method for model reduction. In this section, we briefly review the square root procedure proposed in \cite{TombP87} which is more numerically efficient than its original version. As other projection-based methods, a~balancing projection for system \eqref{eq:psys} must be constructed. This projection helps to balance the controllability and observability energies on each state so that one can easily decide which state component should be truncated. 
 To this end, one has to solve a pair of the generalized Lyapunov equations
\begin{align}
\label{eq:LyapContr}
E(\mu)P(\mu)A^T\!(\mu) + A(\mu)P(\mu)E^T\!(\mu) & = -B(\mu)B^T\!(\mu), \\
\label{eq:LyapObser}
E^T(\mu)Q(\mu)A(\mu) + A^T(\mu)Q(\mu)E(\mu) & = -C^T(\mu)C(\mu),
\end{align}
for the \textit{controllability Gramian} $P(\mu)$ and the \textit{observability Gramian} $Q(\mu)$. In practice, 
these Gramians are computed in the factorized form 
$$
P(\mu) = X(\mu)X^T(\mu), \qquad Q(\mu) = Y(\mu)Y^T(\mu)
$$ 
with $X(\mu)\in\mathbb{R}^{n\times k_c}$ and $Y(\mu)\in\mathbb{R}^{n\times k_o}$. One can show that the eigenvalues of the matrix $P(\mu)E^T(\mu)Q(\mu)E(\mu)$ are real and non-negative \cite{Anto05}. The positive square roots of the eigenvalues of this matrix, $\sigma_1(\mu) \geq \cdots \geq \sigma_{n}(\mu)\geq 0$, are called the \textit{Hankel singular values} of system \eqref{eq:psys}. They can also be determined from
the singular value decomposition (SVD) 
\begin{equation}
\label{eq:SVD}
Y^T(\mu)E(\mu)X(\mu) = [U_1(\mu)\enskip U_0(\mu)]\begin{bmatrix}
\Sigma_1(\mu) & 0\\ 0 & \Sigma_0(\mu)
\end{bmatrix} [V_1(\mu)\enskip V_0(\mu)]^T,
\end{equation}
where $[U_1(\mu)\enskip U_0(\mu)]$ and $[V_1(\mu)\enskip V_0(\mu)]$ are orthogonal, and
$$
\Sigma_1(\mu) = \mbox{diag}(\sigma_1(\mu),\ldots,\sigma_r(\mu)), \quad 
\Sigma_0(\mu) = \mbox{diag}(\sigma_{r+1}(\mu),\ldots,\sigma_{k_{co}}(\mu))
$$
with $k_{co}=\min(k_c,k_o)$.
Then the ROM \eqref{eq:redpsys} is computed by projection 
\begin{equation}\label{eq:RedMatr}
\begin{array}{ll}
\tilde{E}(\mu) = W^T(\mu)E(\mu)T(\mu),& \quad \tilde{A}(\mu) = W^T(\mu)A(\mu)T(\mu), \\
\tilde{B}(\mu) = W^T(\mu)B(\mu), & \quad \tilde{C} = C(\mu)T(\mu),
\end{array}
\end{equation}
where the projection matrices are given by
\begin{equation}
W(\mu) = Y(\mu)U_1(\mu) \Sigma_1^{-1/2}(\mu), \qquad T(\mu) = X(\mu)V_1(\mu)\Sigma_1^{-1/2}(\mu).
\label{eq:WT}
\end{equation}
The $\mathcal{H}_\infty$-error of the approximation is shown to satisfy 
\begin{equation*}
\|H(\cdot,\mu)-\tilde{H}(\cdot,\mu)\|_{\mathcal{H}_\infty} \leq 2\bigl(\sigma_{r+1}(\mu) + \cdots + \sigma_{k_{co}}(\mu)\bigr),
\end{equation*}
where the $\mathcal{H}_\infty$-norm is defined as
\begin{equation*}
\|H\|_{\mathcal{H}_\infty} = \displaystyle{\sup_{\omega\in\mathbb{R}}\|H(\mathrm{i}\omega)\|_2},\mbox{ where } \mathrm{i} = \sqrt{-1},
\end{equation*}
and
\begin{align*}
H(s,\mu)& =C(\mu)(sE(\mu)-A(\mu))^{-1}B(\mu), \\
\tilde{H}(s,\mu)&=\tilde{C}(\mu)(s\tilde{E}(\mu)-\tilde{A}(\mu))^{-1}\tilde{B}(\mu)
\end{align*}
are the transfer functions of systems \eqref{eq:psys} and \eqref{eq:redpsys}, respectively.  

\subsection{Interpolation of Gramians for parametric model order reduction}
Since solving Lyapunov equations is the most expensive step of the balanced truncation procedure, we propose to compute the solution for only a few values of the parameter, and then interpolate those for other values of the parameter. To this end, on the chosen sample grid $\mu_1,\ldots,\mu_q \in \mathcal{D}$, we solve the Lyapunov equations (\ref{eq:LyapContr}) and (\ref{eq:LyapObser}) 
for $P(\mu_i)=P_j=X_j^{}X_j^T$ 
and $Q(\mu_i)=Q_j=Y_j^{}Y_j^T$,
$j=1,\ldots,q$. Note that the ranks of the local Gramians 
$P_j$ and $Q_j$, $j=1,\ldots,q$,
do not need to be the same. Then we define the mappings 
$$
\arraycolsep=2pt
\begin{array}{rclrcl}
P : \mathcal{D} &\to & \R^{n\times n}, &\qquad\qquad Q : \mathcal{D} &\to& \R^{n\times n},\\
 \mu &\mapsto& P(\mu), &\mu &\mapsto& Q(\mu), 
\end{array}
$$
interpolating the data points $(\mu_i, P_i)$ and $(\mu_i, Q_i)$, respectively, as 
\begin{equation*}
P(\mu) = \sum\limits_{j=1}^qw_j(\mu)X_j^{}X_j^T,\qquad Q(\mu) = \sum\limits_{j=1}^qw_j(\mu)Y_j^{}Y_j^T,
\end{equation*}
where $w_j(\mu)$ are some weights that will be detailed in Section \ref{Sec:NumerExam}. To preserve the positive semidefiniteness of the Gramians, we propose to use non-negative weights \cite{Alla03}. 
 This methodology is compatible with the factorization structure since we can write
\begin{align}
P(\mu) &= \sum\limits_{j=1}^q\sqrt{w_j(\mu)}X_{j}\sqrt{w_j(\mu)}X^{T}_{j}\\ \notag
&=\begin{bmatrix}\notag
\sqrt{w_1(\mu)}X_{1}&\cdots & \sqrt{w_q(\mu)}X_{q}
\end{bmatrix}
\begin{bmatrix}  \sqrt{w_1(\mu)}X_{1} &\cdots \sqrt{w_q(\mu)}X_{q}
\end{bmatrix}^T\\ \label{S3GramparaLR}
&= X(\mu)X^T(\mu),
\end{align}
and, similarly,
\begin{equation} \label{S3GramparaLR_2}
	Q(\mu) = Y(\mu)Y^T(\mu),\  \mbox{ where }\ Y(\mu)= \begin{bmatrix}
\sqrt{w_1(\mu)}Y_{1}&\cdots & \sqrt{w_q(\mu)}Y_{q}
\end{bmatrix}.
\end{equation}
Note that the computation of the parametric Gramians is not the ultimate goal. 
 After interpolation, we still have to proceed steps \eqref{eq:SVD} and \eqref{eq:RedMatr} to get the ROM. The computations required by these steps explicitly involve large matrices which may reduce the efficiency of the proposed method. 
To overcome this difficulty, we separate all computations into two stages. The first stage can be expensive but must be independent of $\mu$ so that it can be precomputed.  The second step, where one has to compute the ROM at any new value $\mu\in\mathcal{D}$, 
must be fast. Ideally, its computational complexity should be independent of $n$, the dimension of the initial problem. Such a decomposition is often referred to as an {\em offline-online decomposition} and quite well-known in the reduced basis community \cite{PateR07,HeRS16}.  Details are presented in the next subsection. Before that, we would like to drive the reader's attention to a related work \cite{MassGSSA19}, where 
we considered the problem of interpolating the solution of parametric Lyapunov equations using different interpolation techniques and compared the obtained results. 

\subsection{Offline-online decomposition}
For the offline-online decomposition, we need to add an assumption on the matrices of system (\ref{eq:psys}). We assume here that they can be written as affine combinations of some 
parameter-independent matrices $\{E_i\}_{i = 1, \dots, q_E}$, $\{A_i\}_{i = 1, \dots, q_A}$,  $\{B_i\}_{i = 1, \dots, q_B}$, and $\{C_i\}_{i = 1, \dots, q_C}$ as follows 
\begin{align*}
&E(\mu) = \sum\limits_{i=1}^{q_E}f_i^E(\mu)E_i,\qquad A(\mu) = \sum\limits_{i=1}^{q_A}f_i^A(\mu)A_i,\\
&B(\mu) = \sum\limits_{i=1}^{q_B}f_i^B(\mu)B_i, \qquad C(\mu) = \sum\limits_{i=1}^{q_C}f_i^C(\mu)C_i,
\end{align*}
where $q_E, q_A, q_B, q_C$ are small and the evaluations of $f_i^E,f_i^A,f_i^B,f_i^C$ are cheap. 
Once the interpolated Gramians are available, 
 it follows that
\begin{align}\notag
Y_{}^T(\mu)E(\mu)X_{}(\mu) &= \begin{bmatrix}
\sqrt{w_1(\mu)}Y^{T}_{1}\\ \cdots \\ \sqrt{w_q(\mu)}Y^{T}_{q}
\end{bmatrix} \sum\limits_{i=1}^{q_E}f_i^E(\mu)E_i \begin{bmatrix}
\sqrt{w_1(\mu)}X_{1}&\cdots & \sqrt{w_q(\mu)}X_{q}
\end{bmatrix}\\ \label{S3SVDpara}
&= \sum\limits_{i=1}^{q_E}f_i^E(\mu) 
\begin{bmatrix}
w_{11}(\mu)Y^{T}_{1}E_iX_{1} &\cdots & w_{1q}(\mu)Y^{T}_{1}E_iX_{q}\\
\vdots & \ddots & \vdots \\
w_{q1}(\mu)Y^{T}_{q}E_iX_{1} & \cdots & 
w_{qq}(\mu)Y^{T}_{q}E_iX_{q}
\end{bmatrix},
\end{align}
with $w_{lj}(\mu)=\sqrt{w_l(\mu)w_j(\mu)}$. 
Obviously, all $q_Eq^2$ blocks $Y^{T}_{l} E_i X_{j}$ for $l, j= 1, \ldots, q$ and $i=1, \ldots, q_E$ can be precomputed and stored since they are independent of $\mu$. After computing the SVD of (\ref{S3SVDpara}), the projection matrices 
in \eqref{eq:WT} take the form
\begin{align*}
&W(\mu) = \begin{bmatrix}
\sqrt{w_1(\mu)}Y_{1} & \cdots & \sqrt{w_q(\mu)}Y_{q}
\end{bmatrix}U_1(\mu)\Sigma_1^{-1/2}(\mu),\\
& T(\mu) = \begin{bmatrix}
\sqrt{w_1(\mu)}X_{1} & \cdots & \sqrt{w_k(\mu)}X_{q}
\end{bmatrix} V_1(\mu) \Sigma_1^{-1/2}(\mu). 
\end{align*}
The reduced matrices are then computed as in (\ref{eq:RedMatr}):
\begin{align}\notag
\tilde{E}(\mu) &= W^T(\mu)E(\mu)T(\mu)= \sum\limits_{i=1}^{q_E}f_i^E(\mu)\Sigma_1^{-1/2}(\mu)U_1^T(\mu)\\ \label{S3RedMatPara1}
&\times\begin{bmatrix}
w_{11}(\mu)Y^{T}_{1}E_iX_{1} &\cdots & w_{1q}(\mu)Y^{T}_{1}E_iX_{q}\\
\vdots & \vdots & \vdots \\
w_{q1}(\mu)Y^{T}_{q}E_iX_{1} & \cdots & 
w_{qq}(\mu)Y^{T}_{q}E_iX_{q}
\end{bmatrix}
V_1(\mu) \Sigma_1^{-1/2}(\mu),\\ \notag
\tilde{A}(\mu) &= W^T(\mu)A(\mu)T(\mu)= \sum\limits_{i=1}^{q_A}f_i^A(\mu)\Sigma_1^{-1/2}(\mu)U_1^T(\mu)\\
&\times\begin{bmatrix}
w_{11}(\mu)Y^{T}_{1}A_iX_{1} &\cdots & w_{1q}(\mu)Y^{T}_{1}A_iX_{q}\\
\vdots & \vdots & \vdots \\
w_{q1}(\mu)Y^{T}_{q}A_iX_{1} & \cdots & 
w_{qq}(\mu)Y^{T}_{q}A_iX_{q}
\end{bmatrix}V_1(\mu) \Sigma_1^{-1/2}(\mu),\\ 
\tilde{B}(\mu) &= W^T(\mu)B(\mu)= \sum\limits_{i=1}^{q_B}f_i^B(\mu)\Sigma_1^{-1/2}(\mu)U_1^T(\mu)
\begin{bmatrix}
\sqrt{w_1(\mu)}Y^{T}_{1}B_i\\
\vdots \\
\sqrt{w_q(\mu)}Y^{T}_{q}B_i 
\end{bmatrix},\\
\tilde{C}(\mu) &= C(\mu)T(\mu) \notag \\
&= \sum\limits_{i=1}^{q_C}f_i^C(\mu)\begin{bmatrix}
\sqrt{w_1(\mu)}C_iX_{1} & \cdots & \sqrt{w_q(\mu)}C_iX_{q}\end{bmatrix}V_1(\mu)\Sigma_1^{-1/2}(\mu).\label{S3RedMatPara2}
\end{align}
Again, all matrix blocks, that are independent of $\mu$, can be computed and stored before hand. The offline-online procedure can thus be summarized as follows.
\begin{description}
	\item[\textbf{Offline}] For $\mu_1,\ldots,\mu_q\in\mathcal{D}$,
	\begin{itemize}
		\item solve the Lyapunov equations \eqref{eq:LyapContr} and \eqref{eq:LyapObser} for 
		$P_j\approx X_{j}^{}X_j^T$ and $Q_j\approx Y_j^{}Y_{j}^T$, $j = 1,\ldots,q$;
		\item compute and store all the parameter-independent matrix blocks mentioned in (\ref{S3SVDpara})-(\ref{S3RedMatPara2}).
	\end{itemize}
	\item[\textbf{Online}] Given $\mu\in\mathcal{D}$, 
	\begin{itemize}
		\item assemble precomputed matrix blocks and compute the SVD of (\ref{S3SVDpara});
		\item assemble precomputed matrix blocks and compute the reduced matrices (\ref{S3RedMatPara1})-(\ref{S3RedMatPara2}). 
	\end{itemize}
\end{description}


\section{Interpolation on the manifold $\mathcal{S}_+(k,n)$}\label{Sec:Manifold}

As we have seen above, the low-rank factors of interpolated Grammians obtained by the algebraic approach in \eqref{S3GramparaLR} and \eqref{S3GramparaLR_2} in general have a considerably higher rank than the approximated local Gramians $P_j$ and $Q_j, j = 1,\ldots,q$ obtained by approximately solving the Lyapunov equations. This fact somewhat puts more computational burden on the later steps in model reduction procedure, especially when the number of grid points is large. In many cases when the Gramians do not change so much with respect to the variation of the parameter, we can assume that all approximated local Gramians $P_j, j = 1,\ldots,q$, have a fixed rank. That is, with some relaxation, we can assume that $P_j \in \mathcal{S}_+(k_P,n)$ for $j = 1,\ldots,q$ and $Q_j \in \mathcal{S}_+(k_Q,n)$ for $j = 1,\ldots,q$, where $\mathcal{S}_+(k,n)$ is the set of $n\times n$ positive semidefinite matrices of rank $k$. This set admits a manifold structure \cite{JouBAS10, MassA18}, and therefore our second interpolation method relies on this geometric property.

Informally speaking, a $d$-dimensional manifold is a set $\calM$ that can be mapped locally through a set of bijections, called \emph{charts}, to (an~open subset of) the Euclidean space $\R^d$. Under some additional compatibility assumptions, the collection of charts forms a differentiable structure and the set $\calM$ endowed with this structure is called a $d$-dimensional manifold. The set of charts allows rewriting locally any problem defined on $\calM$ into a problem defined on a subset of $\R^d$. We will see that $\psd{k}{n}$ is a \emph{matrix manifold}, i.e., a manifold whose points can be represented by matrices.

Many matrix manifolds are either \emph{embedded submanifolds} of $\R^{m \times n}$, i.e., the manifold is a subset of the Euclidean space $\R^{m \times n}$, or \emph{quotient manifolds} of $\R^{m \times n}$, where each point of a quotient manifold is a set of equivalent points of $\R^{m \times n}$ for a given equivalence relationship. In each case, the differential structure is inherited from the differential structure on $\R^{m \times n}$. 

As charts are defined locally, they are not very practical for numerical computations. Their use can be avoided by resorting to other tools specific for working on manifolds. The most important for this work are \emph{tangent spaces}, \emph{exponential} and \emph{logarithmic maps}. The tangent space $\calT_x \calM$ is a first order approximation to the manifold $\calM$ around $x \in \calM$, where the point $x$ is called the \emph{foot} of the tangent space. When the tangent spaces are endowed with a \emph{Riemannian metric} (an~inner product) $g_x : \calT_x \calM \times \calT_x \calM \to \R$ smoothly varying with $x$), the manifold is called a \emph{Riemannian manifold}.

The Riemannian metric allows defining geodesics (curves with zero acceleration) on the manifold. This in turn leads to the exponential map which allows mapping tangent vectors to the manifold by following the geodesic starting at the foot of the tangent vector, and whose initial velocity is given by the tangent vector itself. Its reciprocal map is the logarithm map mapping points from the manifold to a given tangent space. For further details on Riemannian manifolds, we refer to \cite{doCarmo1992,Lee2018}. 

\subsection{A quotient geometry of $\mathcal{S}_+(k,n)$}

The manifold $\psd{k}{n}$ is here seen as a quotient manifold $\R_*^{n \times k}/\calO_k$, where $\R_*^{n \times k}$ is the set of full-rank $n \times k$ matrices and $\calO_k$ is the orthogonal group in dimension $k$. This geometry has been developed in~\cite{JouBAS10, MassA18,MassAH19} and has already been used in, e.g.,~\cite{GouMMAJHM17, SabAVGE2019, SzcDBDPM2019} for solving different fitting problems. It relies on the fact that any matrix $A \in \psd{k}{n}$ can be factorized as $A  = YY^{T}$ with $Y \in \R_*^{n \times k}$. As the factorization is not unique, this leads to the equivalence relationship:
\[ Y_1 \sim Y_2  \qquad \text{if and only if} \qquad Y_1 =  Y_2 Q, \ Q \in \calO_k.  \]
For any $Y \in \R_*^{n \times k}$, the set 
\[  [Y] := \{  YQ \enskip :\enskip| Q \in \calO_k\}  \]
of points equivalent to $Y$ 
is called the \emph{equivalence class} of $Y$. The quotient manifold $\R_*^{n \times k}/\calO_k$ is the set of all equivalence classes.

The fact that on the manifold $\R_*^{n \times k}/\calO_k$, any point is a set of points in $\R_*^{n \times k}$ makes it difficult to perform computations directly on elements of $\R_*^{n \times k}/\calO_k$. Instead of manipulating sets of points, most algorithms on quotient manifolds are only manipulating representatives of the equivalence classes.

The tangent space $\calT_Y \R_*^{n \times k}$ to the manifold $\R_*^{n \times k}$ at some point $Y$ is the direct sum of two subspaces: the vertical space $\calV_{Y}$ which is, by definition, tangent to $[Y]$, and the horizontal space $\calH_Y$ which is its orthogonal complement with respect to the Euclidean metric in $\R_*^{n \times k}$. Horizontal vectors will allow  representing tangent vectors to the quotient manifold $\R_*^{n \times k}/ \calO_k$ in a ``tangible way'', i.e., in a way suitable for numerical computations. Indeed, any tangent vector $\xi_{[Y]} \in \calT_{[Y]} \R_*^{n \times k}/ \calO_k$ can be identified to a unique horizontal vector $\bar \xi_Y \in \calH_Y$ in the sense that the two vectors act identically as differentiable operators, see~\cite[\S 3.5.8]{Absil2008}. This vector is called the \emph{horizontal lift} of $\xi_{[Y]}$.

The Riemannian metric is naturally inherited from the Euclidean metric in $\R_*^{n \times k}$ (see~\cite{MassA18}). When defined, the associated exponential map can be written as 
\begin{equation}
\mathrm{Exp}_{[Y]} (\xi_{[Y]}) = [Y + \bar \xi_Y],   \label{eq:exp}
\end{equation}
where  $Y$ is an arbitrary element of the equivalence class $[Y]$, $ \bar \xi_Y$ the unique horizontal lift of the tangent vector $\xi_{[Y]}$. Accordingly, for 
$[Y_1], [Y_2] \in \R_*^{n \times k}/ \calO_k$, 
the logarithm of $[Y_2]$ at $[Y_1]$, $\mathrm{Log}_{[Y_1]}([Y_2])$,  is a vector in 
$\calT_{[Y_1]} \R_*^{n \times k}/ \calO_k$ whose horizontal lift at $Y_1$ is given by
\begin{equation}
\overline{\mathrm{Log}_{[Y_1]}([Y_2])}_{Y_1} = Y_2 Q^{T} - Y_1, \label{eq:log}
\end{equation}
where 
 $Q^{T}$ is the transpose of the orthogonal factor of the polar decomposition of $Y_1^{T} Y_2$, when unique.  We refer the interested reader to~\cite{MassA18} for more information on the domain of definition of these mappings. In all the data sets considered here, we have never faced issues related to ill-definitions of these tools.

\subsection{Curve and surface interpolation on manifolds}

To interpolate the matrices $P_i$ and $Q_i$ on $\psd{k_P}{n}$ and $\psd{k_Q}{n}$, respectively, we consider an intrinsic interpolation technique, general to any Riemannian manifold $\M$. Here, we briefly review it 
in a more 
general framework and refer to~\cite{Gousenbourger2018} and \cite{Absil2016a} for the detailed discussion of curve fitting and surface (i.e., bidimensional) interpolation, respectively.

Consider a 
 Riemannian manifold $\M$ (for instance,
the set of positive semidefinite matrices of size $n$ and rank $k$), and 
a set of data points $d_0\ldots,d_q \in \M$ (like the matrices $P_i$) 
associated to parameter values $t_0, \ldots, t_q \in \R^{\ell}$
(in this case, the values $\mu_i$), where $\ell \in \{ 1, 2\}$ depending on whether one seeks an 
interpolating \emph{curve} or an interpolating \emph{surface}. 

\textbf{Curves.} 
Curve interpolation on $\M$ is often done by encapsulating the interpolation 
into an optimization problem, e.g., one seeks the curve $\bspline:\mathbb{R}\to\M$ minimizing
\begin{equation}  
\label{eq:Eb}
	\min_{\bspline} 
	   \int_{t_{0}}^{t_{q}} \left\|\frac{\mathrm{D}^2 \bspline(t)}{\mathrm{d}t^2}
		\right\|_{\bspline(t)}^2 \mathrm{d}t \ \text{ such that } \ \bspline(t_i) = d_i,
    \quad i = 0,\dots,q,
\end{equation}
where 
the operator $\mathrm{D}^2/\mathrm{d}t^2$ is the Levi-Civita 
covariant derivative of the manifold-valued function $\bspline$ and $\| \cdot \|_{\bspline(t)}$ the norm on the tangent spaces whose foot is along $\bspline(t)$~\cite{Absil2008}.
Different techniques exist to solve this problem, but nearly none of them
tackle~\eqref{eq:Eb} directly on $\M$, as the computational effort 
would be so high that it would not bring any advantages to most 
of the applications.

An efficient way to approximate 
the optimal solution is to transfer the interpolation problem to a 
carefully chosen tangent space $T_x \M$ at a point $x \in \M$, such that 
$T_x \M$ approximates $\M$ in the area where the data points are defined.
The transfer to $T_x \M$ is usually done by mapping the data points to
$T_x \M$ via the logarithmic map or an accurate approximation of it.
As the tangent space is a Euclidean space, solving the Euclidean version
of~\eqref{eq:Eb} is easy and computationally tractable since the Levi-Civita derivative reduces to a 
classical second derivative.  Actually,
the solution can even be written in closed form as it is the natural 
smoothing spline~\cite{Green1993}.
When an approximated curve is computed on $T_x \M$, it is mapped back to
$\M$ via the exponential map or an appropriate retraction, see~\cite{Absil2008,doCarmo1992}
for a detailed exposition on Riemannian geometry.
Curves obtained in this way are noted $\bspline^{TS}(t)$, where the superscript
$TS$ comes from \textbf{T}angent \textbf{S}pace.

It should, however, be noted that the tangent space $T_x \M$ is a good approximation of $\M$ 
only in a close neighborhood of $x$, and in most of the cases, the data 
points cannot all lie in this neighborhood. This is why the so-called
\emph{blended curve} exploits multiple tangent spaces. It is built as a
$C^1$-composite curve~\cite{Gousenbourger2018}
$$\arraycolsep=2pt
\begin{array}{rcl}
  \bspline: [0,1] &\to& \M \\
	t\enskip &\mapsto &  f_i(t-t_i), \quad  \text{when\; } t\in (t_i,t_{i+1}), \;i = 0,\dots,q-1.
	\end{array}
$$
Here,
 $f_i(t)$ is the weighted mean of two curves $\bspline^{TS}(t)$ 
computed, respectively, on the tangent spaces based at $x_i$ and $x_{i+1}$.
This weighted mean is what gives its name (blended) to the technique.

\textbf{Surfaces.}
Interpolation via surfaces is a little bit more intricate. In this work,
we rely on B\'ezier surfaces presented in~\cite{Absil2016a} as a generalization
of Euclidean B\'ezier surfaces~\cite{Farin2002} inspired from the 
manifold generalization of curves to manifolds already presented by Popiel
\emph{et al.}~\cite{Popiel2007}.

Consider a Euclidean space $\R^r$.
A Euclidean B\'ezier curve and B\'ezier surface of degree $K \in \mathbb{N}$ 
are functions $\beta_{K,\ell}$, $\ell=1,2$, 
defined as
\begin{eqnarray}
  \label{eq:euclCurve}
  \beta_{K,1}(\cdot; b_0,\dots,b_K):
    [0,1] &\to&  \R^r \\
		t\enskip &\mapsto & \displaystyle{\sum_{i = 0}^K b_{i} B_{iK}(t)},\nonumber\\
  \label{eq:euclSurface}
  \beta_{K,2}(\cdot,\cdot	; (b_{ij})_{i,j = 0,\ldots,K}) :
    [0,1]^2 & \to & \R^r  \\
				(t^{(1)},t^{(2)}) &\mapsto& \displaystyle{\sum_{i,j = 0}^K b_{ij} B_{iK}(t^{(1)}) B_{jK}(t^{(2)}),} \nonumber
\end{eqnarray}

where $B_{jK}(t) = \binom{K}{j}t^j (1-t)^{K-j}$ are the Bernstein polynomials,
and $b_i \in \R^r$ (resp. $b_{ij} \in \R^r$) are the \emph{control points} of 
the curve (resp. surface).
Since the Bernstein polynomials form a partition of unity, the surface 
can be seen as a convex combination of the control points. Hence, 
equation~\eqref{eq:euclSurface} is equivalent to computing two  
B\'ezier curves~\eqref{eq:euclCurve}: the first one in the $t^{(1)}$ direction, and then in the
$t^{(2)}$ direction, i.e.,
\begin{align*}
  \beta_{K,2}(t^{(1)},t^{(2)}; (b_{ij})_{i,j = 0,\ldots,K})
  & =
  \sum_{i = 0}^K \left( 
    \sum_{i = 0}^K b_{ij} B_{iK}(t^{(1)}) \right) 
    B_{jK}(t^{(2)}) \\
  & = 
  \beta_{K,1}
    \left(
      t^{(2)}; (\beta_{K,1}(t^{(1)}; (b_{ij})_{i=0,\dots,K})_{j=0,\dots,K}
    \right).
\end{align*}
This equivalence permits to easily generalize B\'ezier surfaces to a manifold $\M$
by using the generalization of B\'ezier \emph{curves} 
based on the De Casteljau algorithm, see~\cite{Popiel2007}
for details.

To perform interpolation of data points $d_{ij} \in \M$ associated with
parameter values 
$t_{ij}=[x_i, y_j]^T\in\R^{2}$ one seeks the $C^1$-composite surface 
$$
\arraycolsep=2pt
\begin{array}{rcl}
  \bspline: 
    [x_1,\dots,x_{q_1}]\times[y_1,\dots,y_{q_2}] & \to & \M \\
		(t^{(1)},t^{(2)}) &\mapsto &
    \beta_{K,2}(t^{(1)} - x_k, t^{(2)} - y_l; (b^{kl}_{ij})_{i,j = 0,\dots,K})
		\end{array}
$$
when $t^{(1)} \in (x_k, x_{k+1})$ and $t^{(2)} \in (y_l,y_{l+1})$. Here, 
$\beta_{K,2}(\cdot,\cdot; (b^{kl}_{ij})_{i,j = 0,\dots,K})$ denotes
a~B\'ezier surface on $\M$, and $b^{kl}_{ij}\in\M$ are the control points to be determined such that  
the mean squared second derivative of the piecewise surface is minimized.
This is done with a technique close to the one 
used for curves, i.e., transferring the optimization problem on carefully 
chosen tangent spaces. The only difference here is that the curve itself 
is not computed on the tangent space; instead, the optimality conditions 
obtained on a Euclidean space are generalized to manifolds.
We refer to~\cite{Absil2016a} for a~detailed presentation
of the optimization of the control points, and to \cite{Absil2016} for 
a~complete discussion on the $C^1$-conditions to patch several B\'ezier 
surfaces together.

\section{Numerical examples}\label{Sec:NumerExam}
In this section, we consider  two numerical examples. Before going into detail, we would like to discuss 
the general setting. First, for 
the choice of positive weight coefficients used in the algebraic approach, options are the weights based on distance from the test point to training points and linear splines. Our tests revealed that the latter delivers a smaller error. Moreover, since we have to gather all data to make a big matrix in this method, see \eqref{S3GramparaLR}, too much data may result in inefficiency. Therefore, in the numerical tests, we only use linear splines. In this local interpolation approach, instead of $q$ matrix blocks in each factor of \eqref{S3GramparaLR}, we have only two (resp., four) of them for models with one (resp., two) parameter(s) regardless of the number of training points. One advantage of this approach is that if we want more accuracy by increasing the number of training points, more computation will be required in the offline stage but this makes no changes in the online stage. Thus, this local interpolation is much less affected by the so-called curse of dimensionality when the number of parameters increases compared to the conventional approach. Furthermore, based on the numerical comparisons performed in \cite{MassGSSA19}, in the geometric approach,  we choose the blended curves interpolation technique for the case of one parameter. When the model has two parameters, we use  piecewise B\'ezier surface interpolation. 

To verify the accuracy of ROMs, we compute an approximate $\mathcal{H}_\infty$-norm of the absolute errors in the frequency response defined as
\begin{equation}
\arraycolsep=2pt
\begin{array}{rcl}
\|H(\cdot,\mu)-\tilde{H}(\cdot,\mu)\|_{\mathcal{H}_\infty} & = &
\displaystyle{\sup_{\omega\in\mathbb{R}}\|H(\mathrm{i}\omega,\mu)-\tilde{H}(\mathrm{i}\omega,\mu)\|_2} \\
& \approx &  
\displaystyle{\sup_{\omega_j\in[\omega_{\min},\omega_{\max}]}\|H(\mathrm{i}\omega_j,\mu)-\tilde{H}(\mathrm{i}\omega_j,\mu)\|_2,}
\end{array}
\label{eq:errH}
\end{equation}
where 
$H(s,\mu)$ and $\tilde{H}(s,\mu)$ are the transfer functions of the FOM \eqref{eq:psys} and the ROM \eqref{eq:redpsys}.

For the reference of efficiency, all computations are performed with
MATLAB R2018a on a standard desktop using 64-bit OS
Windows 10, equipped with 3.20 GHz 16 GB Intel Core
i7-8700U CPU.
\subsection{A model for heat conduction in solid material}
This model is adapted from the one used in \cite{KrePT14}.  Consider the heat equation
\begin{align}
\begin{split}\label{Sec8ExHeatEq}
\frac{\partial \vartheta}{\partial t}-\nabla\cdot(\sigma\nabla \vartheta) = f&\quad \mbox{in}\quad \Omega\times (0,T),\\
\vartheta=0&\quad \mbox{on}\; \partial\Omega \times (0,T),
\end{split}
\end{align}
with the heat conductivity coefficient
\begin{equation}
\sigma(\xi) = \left\{ 
\begin{array}{lcl}
1+\mu^{(i)}&\mbox{for}&\xi\in D_i, \; i=1,2,\\
1&\mbox{for}&\xi\in \Omega\backslash (D_1\cup D_2 ),\\
\end{array}
\right.
\end{equation}
where $D_i\subset \Omega=(0,4)^2$, $i=1,2$, are two discs of radius $0.5$ centered at $(1,1)$ 
 and $(3,3)$, respectively, and the parameter $\mu=[\mu^{(1)},\mu^{(2)}]^T$ 
varies in $\mathcal{D}=[1,10]\times [4,10]$. Equation (\ref{Sec8ExHeatEq}) with the source term $f\equiv 1$ is 
discretized using the finite element method with piecewise linear basis functions resulting in 
a~system \eqref{eq:psys} of dimension $n=1580$ with the symmetric positive definite mass matrix 
$E(\mu)\equiv E$ and the stiffness matrix
\begin{equation}
A(\mu) = \mu^{(1)}A_1 + \mu^{(2)} A_2 + A_3,
\end{equation}
where $A_1$ and $A_2$ are symmetric negative semidefinite, and $A_3$ is symmetric negative definite. 
The input matrix $B(\mu)\equiv B\in\mathbb{R}^{n}$ originates from the source function~$f$, and the output matrix is given by
$C(\mu)\equiv C=1/n[1,\ldots, 1]\in\mathbb{R}^{1\times n}$. 
The data were provided by the authors of \cite{KrePT14} for which we would like to thank . 

\begin{figure}[th]
	\begin{minipage}{.5\textwidth}
		\includegraphics[width=1\textwidth]{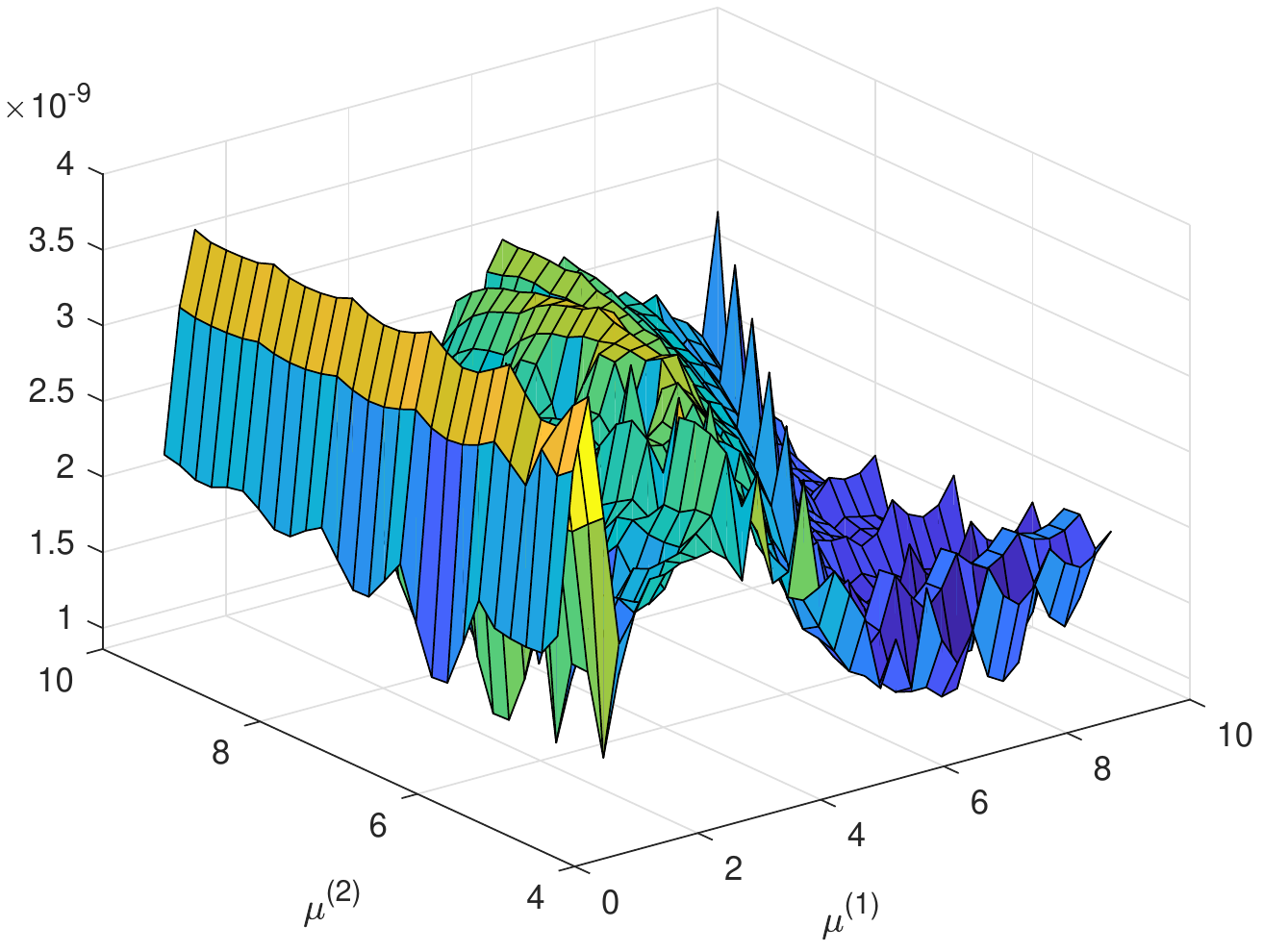}\\
		\includegraphics[width=1\textwidth]{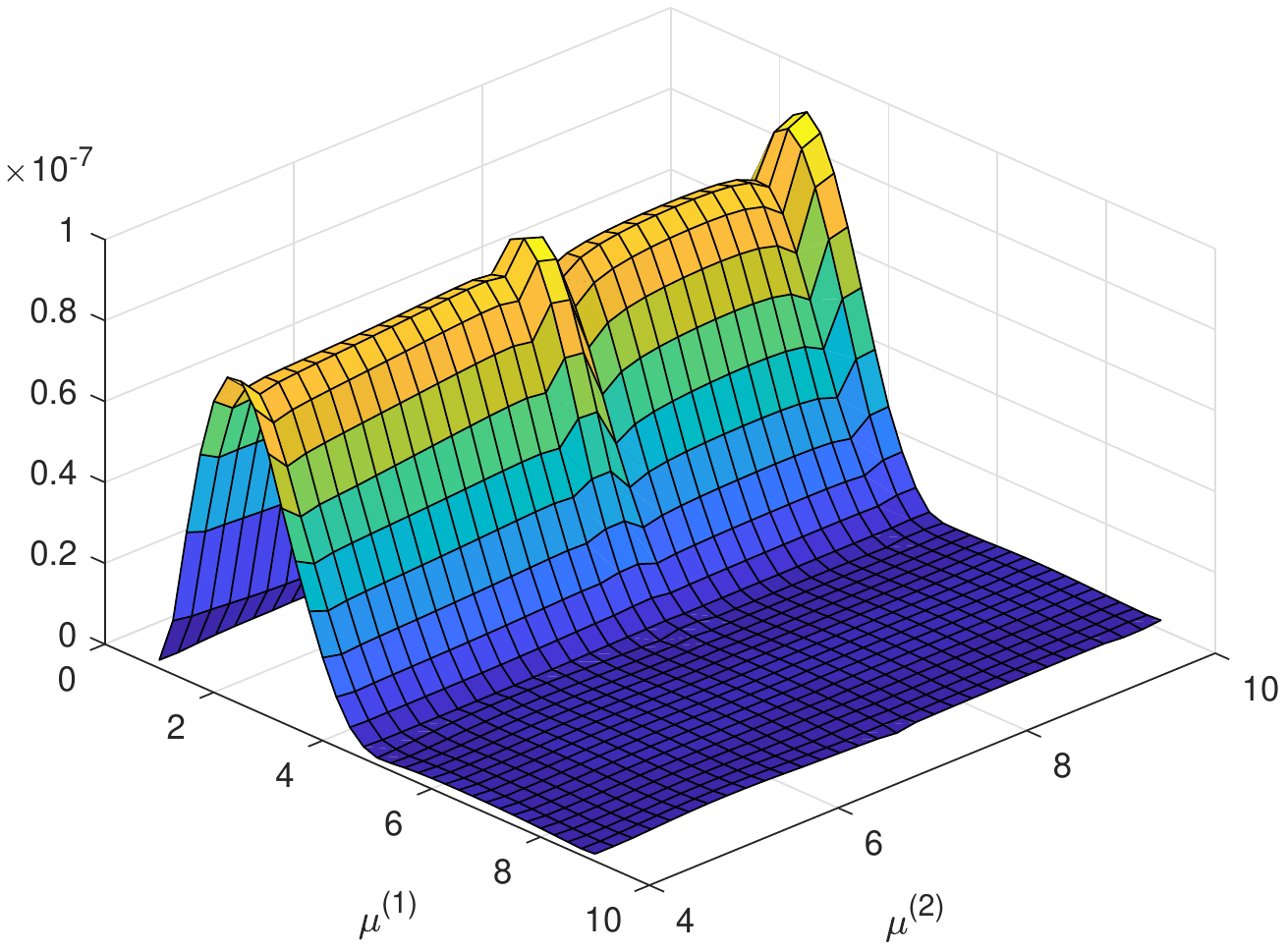}
	\end{minipage}
	\begin{minipage}{.5\textwidth}
		\includegraphics[width=1\textwidth]{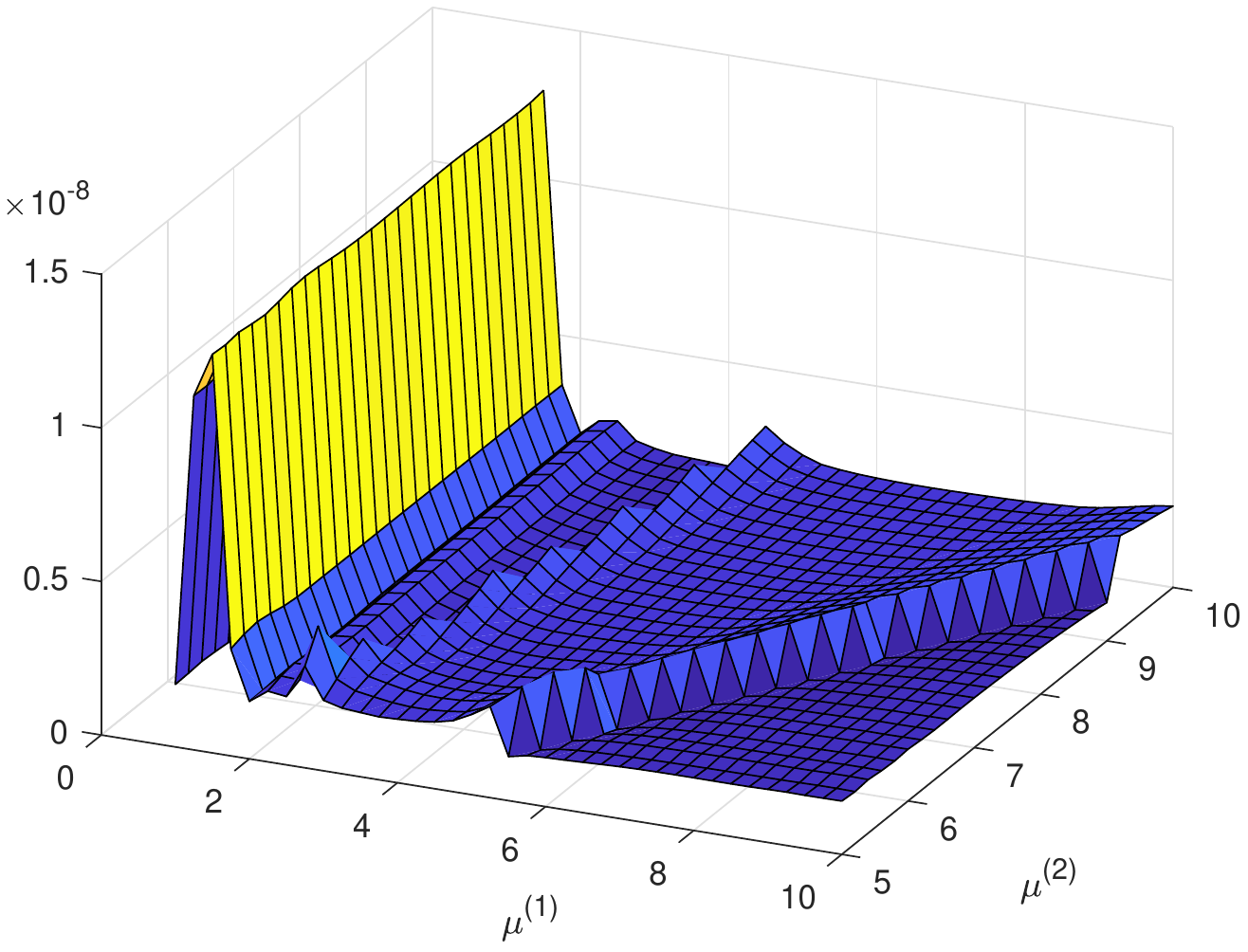} \\
		\includegraphics[width=1\textwidth]{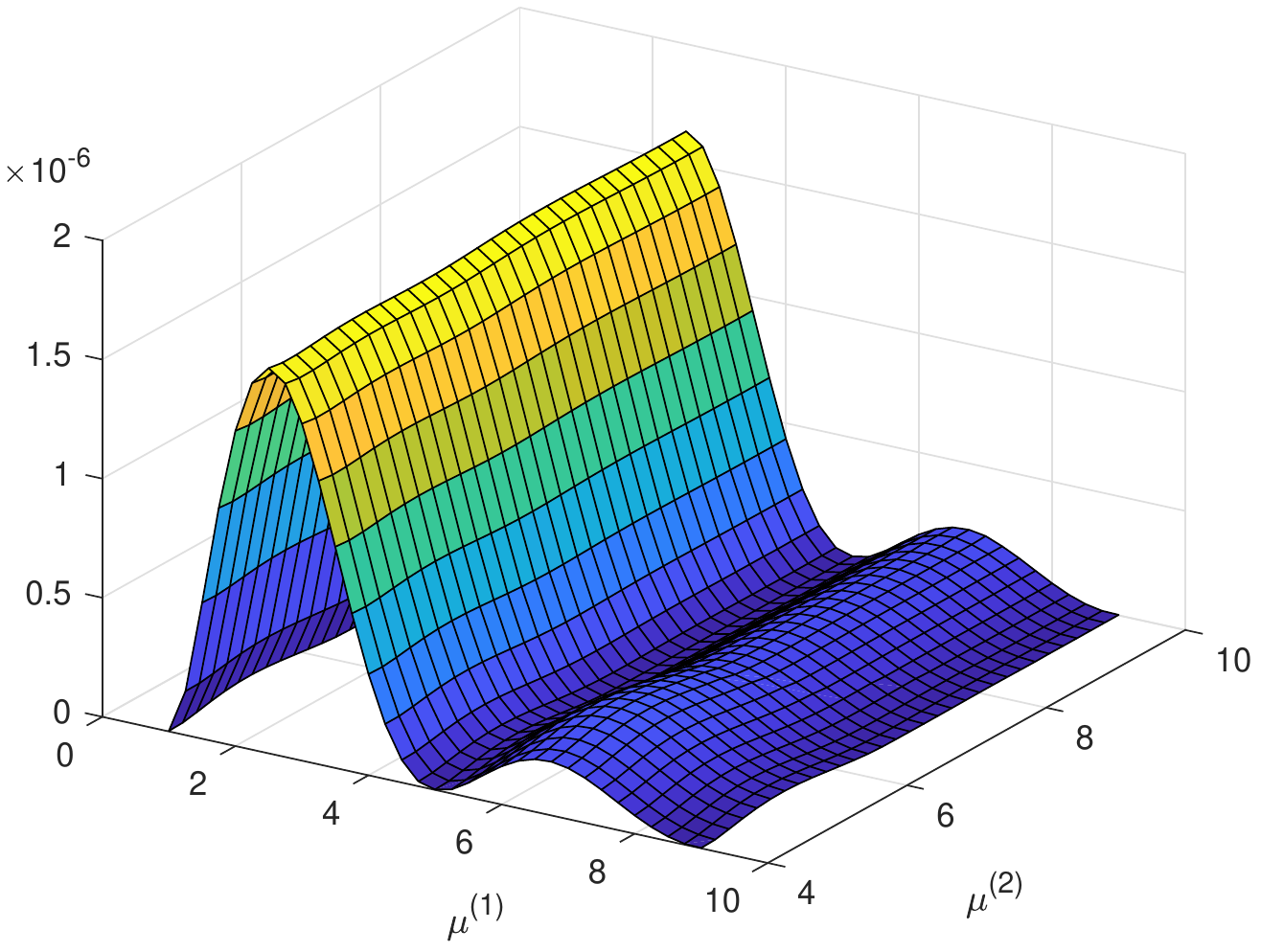}
	\end{minipage}
	\caption{The heat conduction model: 
	absolute errors $\|H(\cdot,\mu)-\tilde{H}(\cdot,\mu)\|_{\mathcal{H}_\infty}$ at test points. 
	Top figures: training grid 
	$[1:1:10]\times[5:1:10]$ and test grid $[1:0.25:10]\times[5:0.2:10]$; 
	bottom figures: training grid 
$[1:4:9]\times[4:3:10]$ and test grid $[1:0.25:9]\times[4:0.2:10]$. 
	The left figures present the errors for the ROMs obtained by the algebraic method and the right figures present that computed by the geometric method.}
	\label{fig:err_example1}
\end{figure}

First, we fix a uniform grid $\mu_1,\dots, \mu_q \in \mathcal{D}$ , which will be specified in the caption of error figures. At those points, we solve \eqref{eq:LyapContr} and \eqref{eq:LyapObser} using the low-rank ADI method \cite{LiW02} with a prescribed tolerance $10^{-10}$. We end up with local approximate solutions whose rank varies from 25 to 27. In order to apply the geometric interpolation method, we truncate them to make all the Gramians of rank 25. In this case, the working manifold is $\mathcal{S}_+(25,1580)$. Note that for the algebraic method presented in Section~\ref{Sec:BT_standard interpolation}, local solutions at training points do not necessarily have the same rank.

\begin{figure}[th]
	\begin{minipage}{.5\textwidth}
	\includegraphics[width=\textwidth]{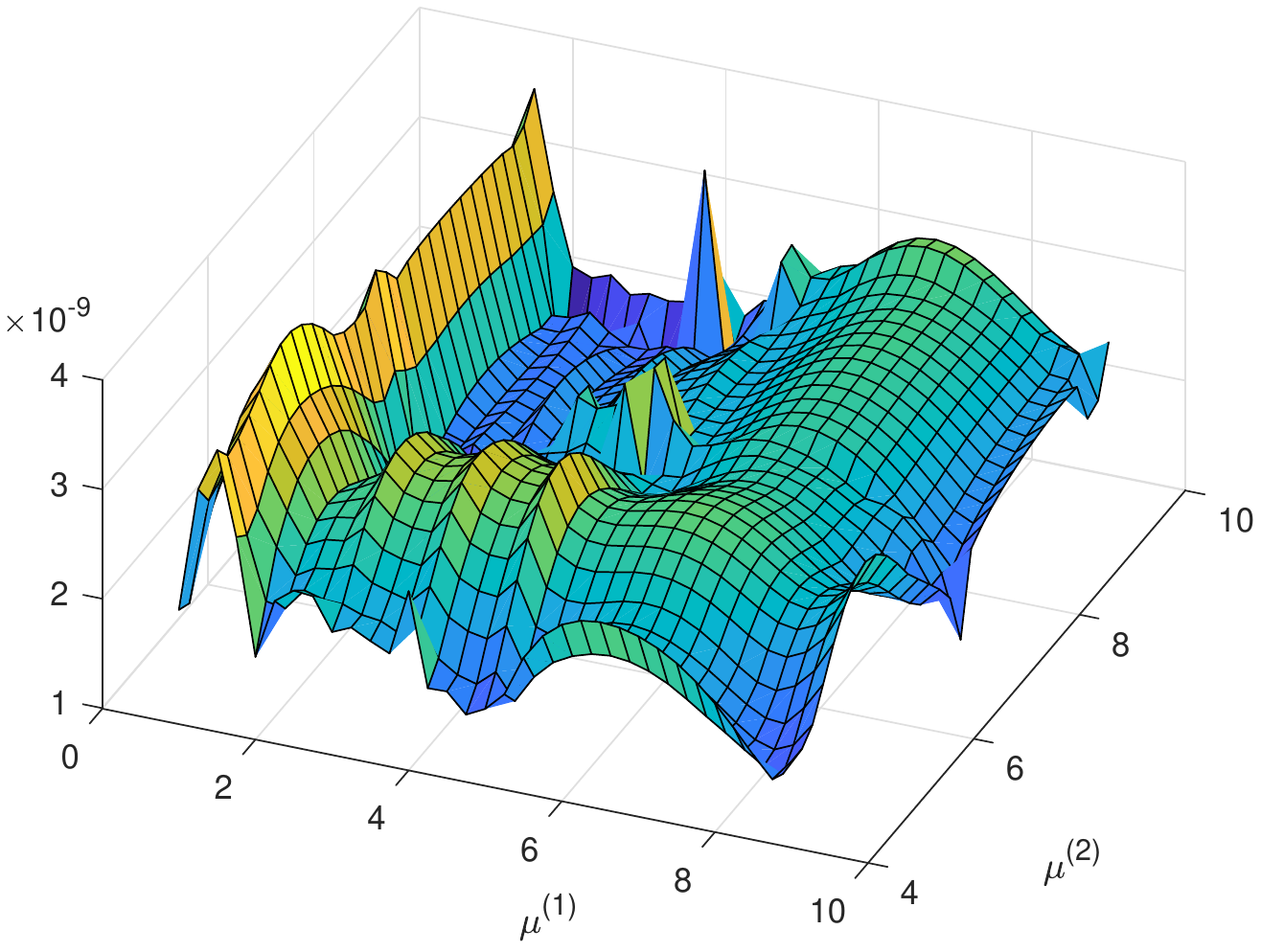}\\
	\end{minipage}
	\begin{minipage}{.5\textwidth}
		\includegraphics[width=\textwidth]{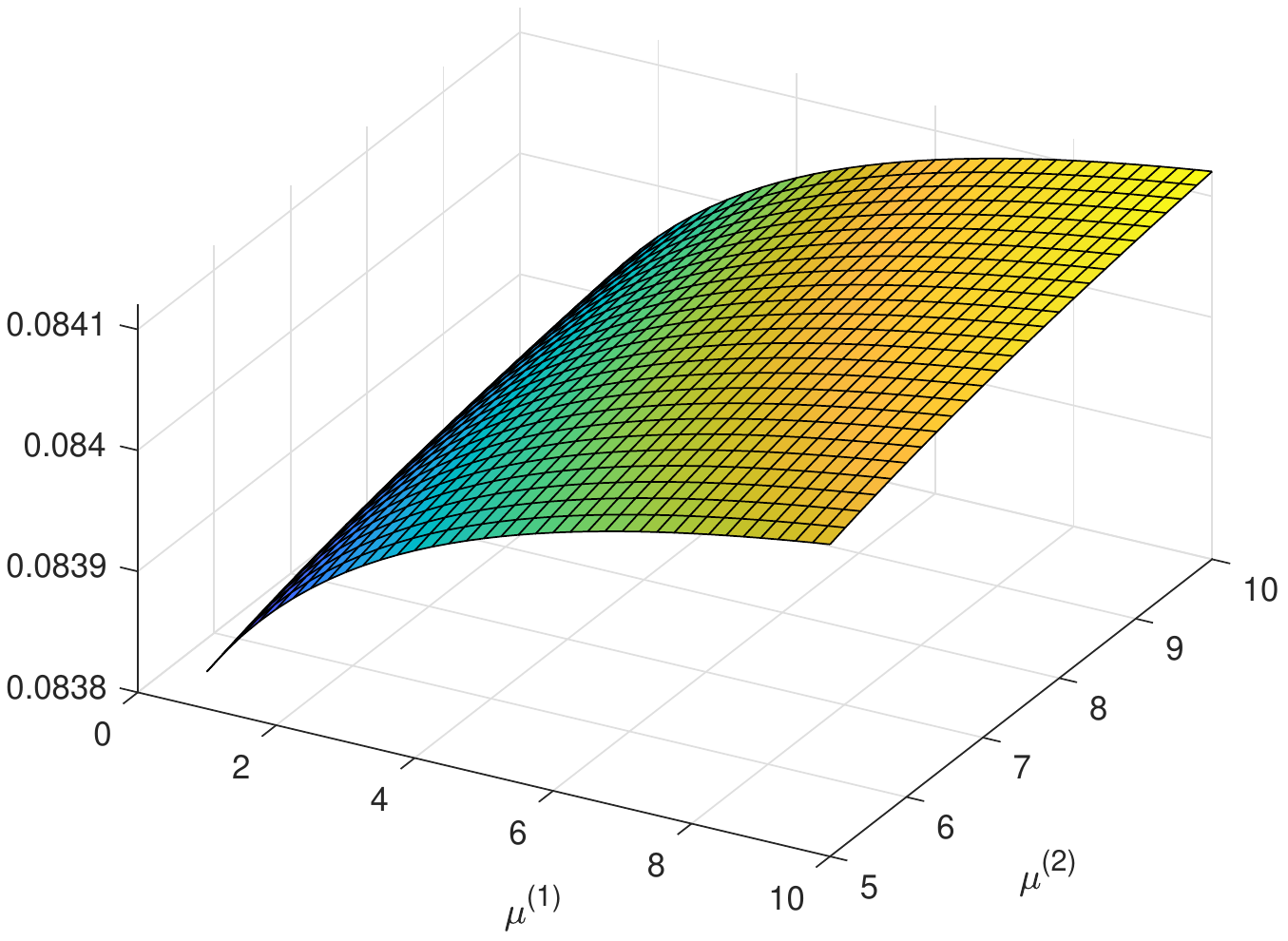} 
	\end{minipage}
	\caption{The heat conduction model: 
	(left) the absolute error $\|H(\cdot,\mu)-\tilde{H}(\cdot,\mu)\|_{\mathcal{H}_\infty}$ with adaptive grid $[1\ 2\ 3\ 4\ 5\ 9]\times[4:3:10]$; (right) the $\mathcal{H}_\infty$-norm of the full-order transfer function on the parameter domain.}
	\label{fig:err2_example1}
\end{figure}

The computed solutions are then employed as the local Gramians to compute the interpolated Gramians which in turn are used to construct the ROM at test points. To determine the reduced order $r$, we use the criterion that $\sigma_r(\mu)/\sigma_1(\mu) < 10^{-8}$ which gives $r$ between 12 to 15 at different test points. 
In Figure~\ref{fig:err_example1}, we plot the approximate absolute errors with respect to $\mathcal{H}_\infty$-norm as defined in \eqref{eq:errH}. For ease of reading numerical results, we simply choose the set of test points as a finer grid of the training grid which will be specified in the caption of the presented 
figures. 
It can be observed that, in the same setting, the algebraic method delivers a slightly smaller error than the geometric one. 
 Moreover, the figures show that the error corresponding to small $\mu_1$ tends to be larger. This suggests that we should use more interpolation data in this area. To this end, we try an adaptively finer grid for the algebraic method and obtain the result as shown in Figure~\ref{fig:err2_example1} (left). Furthermore, to give the reader a view on the relative errors of the method, we plot the $\mathcal{H}_\infty$-norm of the full-order transfer function in Figure~\ref{fig:err2_example1} (right).

We now report the time consumed by the two proposed methods. We will use the second setting that produced the errors as shown in Figure~\ref{fig:err2_example1} (bottom). First, solving two Lyapunov equations at 9 training points needs 2.15 sec. Then, the interpolation of the low-rank solutions of these two equations using the geometric approach at 1023 
test points costs 26.88 sec. From the difference in time consumed, clearly this method can be a good candidate for fast computing the solutions of parametric Lyapunov equations. For model reduction, once the interpolated Gramians 
are available, evaluating ROM at prescribed test points needs 1.47 sec. Meanwhile, for the algebraic approach, the offline stage lasts 0.2 sec and the online one costs 1.33 sec. We summarize these details in Table \ref{tab:TimeConsumedExample11}.



\begin{table}
\caption{The heat conduction model: time consumed by different tasks (sec.)}
\label{tab:TimeConsumedExample11}
	\begin{center}
\begin{tabular}{ p{1cm} p{4cm} | c | c }
	\hline
&	& Geometric approach& Algebraic approach \\\hline
	Offline: & 	
	Solving the Lyapunov equations & & \\ & \quad  at training parameters & 2.15 & 2.15 \\
	&Preparing for interpolation &-  & 0.2 \\
	\hline
	Online:& Interpolation  &26.88 & -\\\
	& Computing the ROMs & 1.47 & 1.33\\
	\hline\noalign{\smallskip}
\end{tabular}
\end{center}
\end{table}


\subsection{An anmometer model}
In the second example, we want to verify the numerical behavior of the proposed methods when applied to fairly large problems. To this end, we consider a~model for a~thermal based flow sensor,  see \cite{MoosRGKH05} and references therein. 
Simulation of this device requires solving a convection-diffusion partial differential equation of the form
\begin{equation}\label{Sec7CDPDEs}
\rho c \frac{\partial \vartheta}{\partial t} = \nabla\cdot(\kappa \nabla \vartheta)-\rho c\mu\nabla \vartheta +\dot{q},
\end{equation}
where $\rho$ denotes the mass density, $c$ 
 is the specific heat, $\kappa$
   is the thermal conductivity, 
$\mu$
 is the fluid velocity, $\vartheta$ is the temperature, and $\dot q$ is the heat flow into the system 
caused by the heater. The considered model is restricted to the case $\rho = 1, c = 1, \kappa = 1$ and $\mu \in [0,1]$ which corresponds to the 1-parameter model. The finite element discretization of 
\eqref{Sec7CDPDEs} leads to system \eqref{eq:psys} of order $n=29008$ with the symmetric positive definite mass matrix $E(\mu) \equiv E$ 
and the stiffness matrix $A(\mu)=A_1+\mu A_2$, where 
$A_1$ is symmetric negative definite, $A_2$ is non-symmetric  negative semidefinite.
The input matrix $B\in\mathbb{R}^n$ and the output matrix $C\in\mathbb{R}^{1\times n}$ are parameter-independent.
The reader is referred to \cite{morwiki_anemom} and references therein for more detailed descriptions and  numerical data.

\begin{figure}[th]
\centering
		\includegraphics[width=0.5\textwidth]{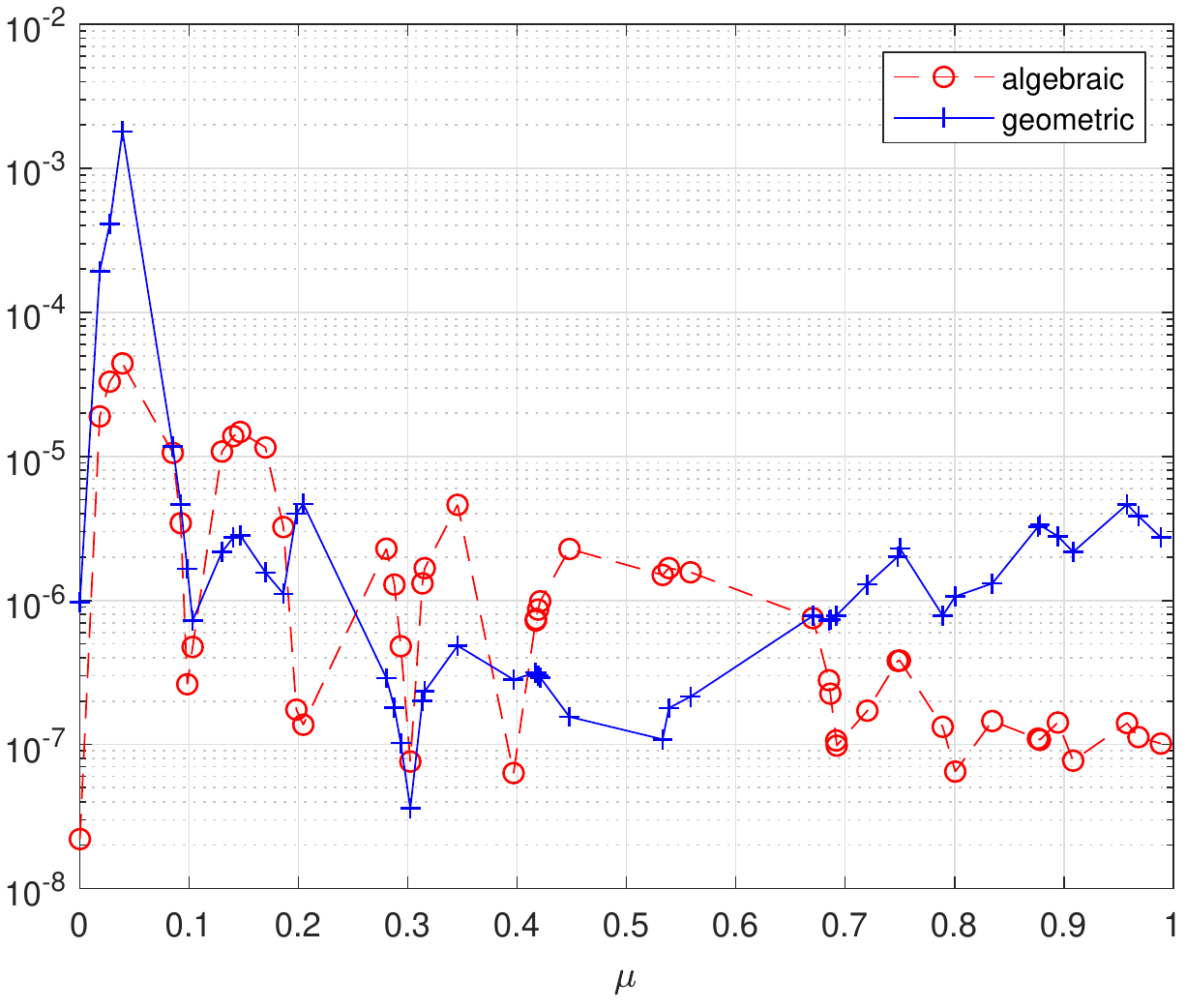}
	\caption{The anemometer model: absolute errors $\|H(\cdot,\mu)-\tilde{H}(\cdot,\mu)\|_{\mathcal{H}_\infty}$ at test points.}
	\label{fig:example2}
\end{figure}

For this model, we use the training grid as $[0:0.1:1]$ while the test grid with 50 points is randomly generated within the range of the parameter domain. The tolerance for the low-rank ADI solver is $10^{-9}$ and that for balancing truncation is $10^{-7}$.
The resulting ROMs have different reduced orders at test points: the ROMs produced by the geometric approach have orders between 9 and 17 while that obtained by the algebraic approach between 16 and 27. The absolute errors are visualized in Figure~\ref{fig:example2}. 
One can see that on some parts of the parameters domain, the geometric approach provides better approximations than that computed by the algebraic methods, while on the others, we observe the reverse results. 

The time consumed by different tasks is summarized in Table~\ref{tab:TimeConsumedExample21}.

\begin{table}
\caption{The heat anemometer model: time consumed by different tasks (sec.)}
\label{tab:TimeConsumedExample21}
	\begin{center}
\begin{tabular}{ p{1cm} p{4cm} | c | c }
	\hline
&	& Geometric approach& Algebraic approach \\\hline
	Offline: & 	
	Solving the Lyapunov equations & & \\ & \quad  at training parameters & 199.18 & 199.18 \\
	&Preparing for interpolation &-  & 25.55 \\
	\hline
	Online:& Interpolation  &18.40 & -\\\
	& Computing the ROMs & 0.81 & 0.06\\
	\hline\noalign{\smallskip}
\end{tabular}
\end{center}
\end{table} 


%

\section{Conclusion}\label{Sec:Concl}
We presented two methods for interpolating the Gramians of parameter-dependent linear dynamical systems for using in parametric balanced truncation model reduction. The first method is merely based on linear algebra which takes no geometric structure of data into account. Thanks to  simplicity, it can be combined with the reduction process which enables an offline-online decomposition. This decomposition in turn accelerates the MOR process  in the online stage which suits very well in parametric settings. Moreover, it is more flexible with the change of parameter values and easier to implement. Meanwhile, the second method exploits the positive semidefiniteness of the data set and recent developments in matrix manifold theory. It reformulates the problem as interpolation on the underlying manifold and relies on the advanced techniques involving interpolating on different tangent spaces and blending the resulting objects to preserve the geometric structure as well as the regularity of data. This method is a good choice for fast interpolating the low-rank solutions of parametric Lyapunov equations and expected to work well if the numerical rank of such solutions does not change much. While the implementation of the geometric approach is challenging, it can result in lower reduced order  
as it often offers  better approximation to the solution of the Lyapunov equations. 


\section*{Acknowledgments}
This work was supported by the Fonds de la Recherche Scientifique - FNRS and the Fonds Wetenschappelijk Onderzoek - Vlaanderen under EOS Project no
30468160.

\backmatter


\end{document}